\documentclass[12pt]{amsart}
\usepackage[alphabetic]{amsrefs}
\usepackage{mathdots, blkarray}
\usepackage{multicol}
\usepackage{graphicx}
\usepackage{amscd}
\usepackage{upgreek}
\usepackage{stmaryrd}
\usepackage{longtable}
\usepackage[T1]{fontenc}
\usepackage{latexsym, amsmath, amssymb, amsthm}
\usepackage[svgnames]{xcolor}
\usepackage{rsfso}
\usepackage[mathscr]{eucal}
\usepackage{mathtools}
\usepackage{mathptmx}
\usepackage{titletoc}
\usepackage{wrapfig}
\usepackage{float}
\usepackage{xypic}
\usepackage{microtype}
\usepackage{dsfont}
\usepackage{xcolor}
\usepackage{color}
\usepackage[colorlinks = true,
            linkcolor  = DarkBlue,
            urlcolor   = DarkRed,
            citecolor  = DarkGreen]{hyperref}
\usepackage{hyperref}
\allowdisplaybreaks	
\usepackage[centering, includeheadfoot, hmargin=1.0in, tmargin=0.5in, 
  bmargin=0.6in, headheight=6pt]{geometry}
\newtheorem{theorem}{Theorem}[section]
\newtheorem{prop}[theorem]{Proposition}
\newtheorem{defn}[theorem]{\rm\textsc{Definition}}
\newtheorem{lem}[theorem]{Lemma}
\newtheorem{coro}[theorem]{Corollary}

\newtheorem{thm}[theorem]{Theorem}

\newtheorem{rem}[theorem]{\rm\textsc{Remark}}
\newtheorem{exam}[theorem]{\rm\textsc{Example}}

\DeclareMathOperator{\Hom}{Hom}

\newcommand{\A}{\mathcal{A}} 
 
\newcommand{\F}{\mathbb{F}} 
\newcommand{\C}{\mathbb{C}} 
 
\newcommand{\Z}{\mathbb{Z}} 
\newcommand{\N}{\mathbb{N}} 
 
\newcommand{\gl}{\mathfrak{gl}} 

\newcommand{\ra}{\longrightarrow}

\newcommand{\ha}{\widehat}  
\newcommand{\hbo}{$\hfill\Diamond$}

\begin{document}
\title{Cohomology of left-symmetric color algebras} 
\def\shorttitle{Cohomology of left-symmetric color algebras}

\author{Yin Chen}
\address{School of Mathematics and Physics, Jinggangshan University,
Ji'an 343009, Jiangxi, China \& Department of Finance and Management Science, University of Saskatchewan, Saskatoon, SK, Canada, S7N 5A7}
\email{yin.chen@usask.ca}

\author{Runxuan Zhang}
\address{Department of Mathematics and Information Technology, Concordia University of Edmonton, Edmonton, AB, Canada, T5B 4E4}
\email{runxuan.zhang@concordia.ab.ca}

\begin{abstract}
We develop a new cohomology theory for finite-dimensional left-symmetric color algebras and their finite-dimensional bimodules, establishing a connection between Lie color cohomology and left-symmetric color cohomology. We prove that the cohomology of
a left-symmetric color algebra $A$ with coefficients in a bimodule $V$ can be computed by a lower degree cohomology of the corresponding Lie color algebra with coefficients in $\Hom(A,V)$, generalizing a result of Dzhumadil'daev in right-symmetric cohomology. We also explore the varieties of two-dimensional and three-dimensional left-symmetric color algebras.
\end{abstract}

\date{\today}
\maketitle \baselineskip=16.5pt

\dottedcontents{section}[1.16cm]{}{1.8em}{5pt}
\dottedcontents{subsection}[2.00cm]{}{2.7em}{5pt}

\section{Introduction}
\setcounter{equation}{0}
\renewcommand{\theequation}
{1.\arabic{equation}}
\setcounter{theorem}{0}
\renewcommand{\thetheorem}
{1.\arabic{theorem}}

\noindent Cohomology of an algebraic structure occupies a central position in understanding representations, extensions, deformations of these algebras, that dates back to many classical work, for example, \cite{Hoc45} for associative algebras, \cite{CE48} and \cite{HS53} for Lie groups and Lie algebras, and \cite{Ger63} for associative rings. 
A cohomology theory of right-symmetric algebras was established in \cite{Dzh99} which reveals a close relationship between
Lie cohomology and right-symmetric cohomology. A theorem of  Dzhumadil'daev  (\cite[Theorem 3.4]{Dzh99}), showing that
the $(n+1)$-th cohomology of a right-symmetric algebra $A$ in coefficients $V$ is isomorphic to the $n$-th cohomology of 
the corresponding Lie algebra $[A]$ in coefficients $\Hom(A,V)$, serves as our motivating example. We develop a cohomology theory for finite-dimensional left-symmetric color algebras and their bimodules in terms of Dzhumadil'daev's philosophy, emphasizing to establish a connection between Lie color cohomology and left-symmetric color cohomology. 

Left-symmetric color algebras, as a generalization of left-symmetric algebras and left-symmetric superalgebras, are
nonassociative graded algebras over a field $k$ graded by an abelian group $G$ with the product satisfying an $\upepsilon$-left-symmetric identity for a skew-symmetric bicharacter $\upepsilon$ of $G$ over $k$. 
Left-symmetric color algebras are Lie admissible color algebras, having a very close relationship with Lie color algebras. 
In contrast to the relatively well-explored theory of Lie color algebras (see for example \cite{Fel01} and \cite{LT23}), the theory of left-symmetric color algebras is quite new and just few results are known; see \cite{CCD17} for constructions of left-symmetric color algebras and \cite{NBN09} for a cohomology and homology theory for left-symmetric color algebras.

The cohomology theory we aim to develop in this article differs from the cohomology theory occurred in \cite{NBN09}. 
For any cohomology theory, it is well-known that cochain spaces and coboundary operators are 
fundamental and essential ingredients. Given a finite-dimensional  left-symmetric color algebra $A$ with respect to $(G,\upepsilon)$ and a finite-dimensional $A$-bimodule $V$, here we define the cochain space $C^n(A,V)$ to be the $G$-graded vector space of all homogeneous linear maps from
the tensor product of $A$ and $n$-th $\upepsilon$-exterior powers of $A$ to $V$, while the cochain spaces in \cite{NBN09} were chosen as the spaces of multilinear maps of $A$ to $V$.
Actually, our choice seems to be more natural as a generalization of the cochain space in the right-symmetric cohomology; see \cite[Section 3.1]{Dzh99}.

Undoubtedly, computations of cohomology spaces and interpretations of low-degree cohomologies are  core topics in a cohomology theory. We define the $0$-th cohomology space $H^0(A,V)$ as the space of the ``invariants'' of $A$ on $V$ and we observe that the first cohomology space can be understood as the space of outer $\upepsilon$-derivations of $A$ to $V$; see Remark \ref{rem3.7} and Example \ref{coh1}. For higher cohomology spaces $H^{n+1}(A,V)$ ($n\in\N^+$), 
the following main result of this article demonstrates that $H^{n+1}(A,V)$  can be computed by 
the lower degree cohomology $H^{n}([A],C^1(A,V))$ of the corresponding Lie color algebra $[A]$ in coefficients $C^1(A,V)$. 

\begin{thm} \label{maint}
Let $A$ be a finite-dimensional left-symmetric color algebra with respect to a skew-symmetric bicharacter $\upepsilon$ of an abelian group $G$ over a field $k$ and let $V$ be a finite-dimensional $A$-bimodule. Then the natural isomorphism $\varphi$ defined in (\ref{ni}) below induces an isomorphism of  $G$-graded vector spaces:
\begin{equation}
\label{ }
H^{n+1}(A,V)\cong H^n([A], C^1(A,V))
\end{equation}
for all $n\in\N^+$.
\end{thm}

The main content of this article is drawn from \cite{Zha13}, which was posted on arXiv in 2013 but  not formally submitted for publication in any journals. Several years later, after the first author joined the research group, we started to re-organize the manuscript and materials, added more examples and details to the proofs, and eventually decided to split the manuscript into the present article that focuses on cohomology theory of left-symmetric color algebras, and the upcoming one \cite{CZ24} regarding the aspect of
deformation theory and applications. We also notice that cohomology and deformation theory of some other generalizations of left-symmetric algebras have been studied recently; see for example, \cite{BHCMN19} for BiHom-left-symmetric algebras, \cite{LMS22} for Hom-left-symmetric algebras, \cite{ZH22} for left-symmetric conformal algebras, and \cite{BA24} for left-symmetric superalgebras.  Some introductory materials about left-symmetric (or pre-Lie) algebras can be found in \cite{Bai20}.

\vspace{2mm}
\noindent \textbf{Layout}. We organize this article as follows.  Section \ref{sec2} contains fundamentals about left-symmetric color algebras, including 
skew-symmetric bicharacters of an abelian groups over a field, Lie color algebras, and bimodules of left-symmetric color algebras. 
We provide lots of examples to illustrate these concepts, and especially, we define a bimodule structure on the tensor product of two bimodules of a 
left-symmetric color algebra (Proposition \ref{prop2.14}) and a right-trivial bimodule structure on the first cochain space (Proposition \ref{bim1}) that plays a crucial role in defining the left-symmetric color cohomology. 

In Section \ref{sec3}, we develop our cohomology theory for a given left-symmetric color algebra $A$ and its bimodule $V$. 
We give a detailed interpretation for the $0$-th or the first spaces of cocycles (coboundaries, cohomologies).  The bulk of this section is devoted to the proof of Theorem \ref{maint}. 

Section \ref{sec4} contains a rough classification of 2-dimensional left-symmetric color algebras and a discussion about  3-dimensional left-symmetric color algebras and their cohomologies.

\vspace{2mm}
\noindent \textbf{Conventions}.  Throughout this article, all algebras and modules are assumed to be finite-dimensional.  We use $|x|$ to denote the degree of a homogeneous element $x$. The symbol $\widehat{y}$ denotes that an element $y$ was deleted. 
Some notations are standard, for example, $k^\times$ denotes the multiplicative group of a field $k$, $\N=\{0,1,2,\dots\}$,
and $\N^+=\{1,2,3,\dots\}$.

\section{Left-Symmetric Color Algebras}\label{sec2}
\setcounter{equation}{0}
\renewcommand{\theequation}
{2.\arabic{equation}}
\setcounter{theorem}{0}
\renewcommand{\thetheorem}
{2.\arabic{theorem}}

\subsection{Skew-symmetric bicharacters of abelian groups} Let $k$ be a field of any characteristic and $G$ be an abelian group. 
We denote by $k^\times$ the multiplicative group of $k$. A map $\upepsilon:G\times G\ra k^\times$
is called a \textit{bicharacter} of $G$ over $k$ if it is biadditive, i.e., $\upepsilon(a,b+c)=\upepsilon(a,b)\cdot \upepsilon(a,c)$ and $\upepsilon(a+b,c)=\upepsilon(a,c)\cdot \upepsilon(b,c)$, for all $a,b,c\in G$. 

A bicharacter $\upepsilon$ of $G$ is said to be \textit{skew-symmetric} if  $\upepsilon(a,b)\cdot\upepsilon (b,a)=1$ for all $a,b\in G$. 
We use $B_k(G)$ to represent the set of all skew-symmetric bicharacters of $G$ over $k$. Clearly,
$\upepsilon(0,a)=1$ and $\upepsilon(a,a)$ is either $1$ or $-1$ for $\upepsilon\in B_k(G)$ and all $a\in G$.

\begin{exam}\label{exam2.1}
{\rm The trivial map on $G\times G$ defined by $\upepsilon(a,b):=1$ for all $a,b\in G$ is clearly a 
skew-symmetric bicharacter of $G$ on $k$. This means that $B_k(G)\neq\emptyset$ for each pair $(G,k)$.
\hbo}\end{exam}

\begin{exam}\label{exam2.2}
{\rm Suppose $G$ denotes the additive group of the prime field $\F_2=\{0,1\}$ of order $2$ and 
the characteristic of $k$ is not $2$. The map 
$\upepsilon:G\times G\ra k^\times$ defined by $$\upepsilon(a,b):=(-1)^{ab}$$ for all $a,b\in G$ is a 
skew-symmetric bicharacter of $G$ over $k$. In fact, 
$\upepsilon(a,b+c)=(-1)^{a(b+c)}=(-1)^{ab}\cdot (-1)^{ac}=\upepsilon(a,b)\cdot \upepsilon(a,c)$, and similarly, we have
$\upepsilon(a+b,c)=\upepsilon(a,c)\cdot \upepsilon(b,c)$. Moreover, $\upepsilon(a,b)\cdot\upepsilon (b,a)=
(-1)^{ab}\cdot (-1)^{ba}=(-1)^{ab}\cdot (-1)^{ab}=(-1)^{2ab}=(-1)^0=1$.
\hbo}\end{exam}

More generally, we have

\begin{exam}\label{exam2.3}
{\rm
Let $p$ be a prime and $\F_p$ be the prime field of order $p$. Suppose $G$ denotes the 
the additive group of the $n$-dimensional vector space $\F_p^n$ over $\F_p$ and 
there exists a nontrivial $p$th root  $\lambda_p$ of unity  in $k$. Given a skew-symmetric bilinear form $\varphi:
\F_p^n\times \F_p^n\ra\F_p$, the map $\upepsilon:\F_p^n\times \F_p^n\ra k^\times$ defined by 
$\upepsilon(a,b):=\lambda_p^{\varphi(a,b)}$ for all $a,b\in \F_p^n$ is a 
skew-symmetric bicharacter of $\F_p^n$ over $k$; see \cite[Section 5]{Sch79} for more details. 

In particular, if $p=2, n=3$, and $k=\C$, then $\lambda_p=-1$. 
\begin{enumerate}
  \item The skew-symmetric bilinear form $\varphi:
\F_2^3\times \F_2^3\ra\F_2$ defined by $\varphi(a,b)=1$ if $a\neq b\in \F_2^3$, and $\varphi(a,a)=0$ for all $a\in\F_2^3$, induces 
a skew-symmetric bicharacter  $\upepsilon:\F_2^3\times \F_2^3\ra\C^\times$  determined by 
$$(\upepsilon(a,b))_{3\times 3}=\begin{pmatrix}
   1   &  -1&-1  \\
    -1  & 1&-1\\
    -1&-1&1 
\end{pmatrix}.$$
  \item The skew-symmetric bilinear form $\varphi:
\F_2^3\times \F_2^3\ra\F_2$ defined by $\varphi(a,b)=0$ if $a\neq b\in \F_2^3$, and $\varphi(a,a)=1$ for all $a\in\F_2^3$, induces 
a skew-symmetric bicharacter  $\upepsilon:\F_2^3\times \F_2^3\ra\C^\times$  determined by 
$$(\upepsilon(a,b))_{3\times 3}=\begin{pmatrix}
   -1   &  1&1  \\
    1  & -1&1\\
    1&1&-1 
\end{pmatrix}.$$
\end{enumerate}
Moreover, all injective skew-symmetric bicharacters of $\F_2^3$ over $\C$ have been classified as a list of $44$ cases in \cite[Section 2]{Sil97}. These two bicharacters above were labelled as $\upepsilon_{18}$ and $\upepsilon_{27}$ respectively in the list. 
\hbo}\end{exam}

\begin{rem}{\rm
In the previous two examples, $(-1)^0$ and $\lambda_p^{0}$ both are well-defined. In fact, assume that we identify the prime field
$\F_p$ with $\Z_p=\{\bar{0},\bar{1},\dots,\overline{p-1}\}$, the ring of residue classes modulo $p$. Then 
$0\in\F_2$ denotes the set of all even integers and thus $(-1)^0$ can be defined as $(-1)^{2r}=1$ for any $r\in\Z$ in 
this case. For general $p$, $\lambda_p^{0}$ can be defined as $(\lambda_p)^{p\cdot r}=1^r=1$ for all $r\in\Z$.
\hbo}\end{rem}

\subsection{Left-symmetric color algebras} 
Given an abelian group $G$, a nonassociative algebra $A$ over a field $k$ is called \textit{$G$-graded} if there exists a decomposition $A=\bigoplus_{a\in G}A_a$ as vector space  such that $A_a\cdot A_b\subseteq A_{a+b}$ for all
$a,b\in G$. We say that an element $x$ of $A$ is \textit{homogeneous of degree $a$} if $x\in A_a$ and denoted as $|x|=a$.
A $G$-graded nonassociative algebra $A=\bigoplus_{a\in G}A_a$ is called a \textit{left-symmetric color algebra} if there exists an $\upepsilon\in B_k(G)$ such that 
\begin{equation}
\label{lsa}
(xy)z-x(yz)=\upepsilon(|x|,|y|)\Big((yx)z-y(xz)\Big)
\end{equation}
for all homogeneous elements $x,y$, and $z\in A$.

A  left-symmetric color algebra is also called a \textit{generalized left-symmetric algebra} in \cite{Zha13}, following
Manfred Scheunert's steps in \cite{Sch79} in which \textit{Lie color algebras} were named as \textit{generalized Lie algebras}. However, seems the name of left-symmetric color algebras has become more popular recently, see for example \cite{CCD17}.

Apparently, left-symmetric color algebras generalize the notion of left-symmetric algebras. In particular, if $A=\bigoplus_{a\in G}A_a$ is a left-symmetric color algebra, then the subspace $A_0$ is a left-symmetric algebra and each subspace $A_a$ is an $A_0$-module; see \cite[Definition 7.8]{Bai20} for modules of left-symmetric algebras. 

The following two examples demonstrates that the class of left-symmetric color algebras not only includes
left-symmetric algebras and $G$-graded left-symmetric algebras but also includes left-symmetric superalgebras as subclasses. 

\begin{exam}\label{}
{\rm Suppose that $\upepsilon$ is the trivial map defined in Example \ref{exam2.1}. Then a left-symmetric color algebra with respect to $(G,\upepsilon)$ is exactly a $G$-graded left-symmetric algebra. Moreover, if the $G$-grading is trivial, then the category of $G$-graded left-symmetric algebras coincides with the category of left-symmetric algebras.
\hbo}\end{exam}

\begin{exam}\label{}
{\rm Let $\upepsilon$ be the skew-symmetric bicharacter defined in Example \ref{exam2.2}.
Then a left-symmetric color algebra with respect to $\upepsilon$ is a left-symmetric superalgebra; see for example
 \cite{DZ23}, \cite{ZHB11}, and \cite{ZB12} for some recent progress in the theory of left-symmetric superalgebras.
\hbo}\end{exam}

\begin{exam}\label{}
{\rm An associative $G$-graded algebra is a left-symmetric color algebra with the same multiplication and any skew-symmetric bicharacter $\upepsilon\in B_k(G)$. See \cite[Example 2.3]{Sch83} for the concept of  associative $G$-graded algebras.
\hbo}\end{exam}

\begin{exam}\label{}
{\rm A linear map $f$ between two $G$-graded vector spaces $V=\bigoplus_{a\in G} V_a$ and 
$W=\bigoplus_{a\in G} W_a$ is called \textit{homogeneous of degree $b\in G$} if $f(V_a)\subseteq W_{a+b}$ for all $a\in G$.

Suppose $G$ is an abelian group such that $V_a=\{0\}$ and $W_{a}=\{0\}$ for all but finitely many $a\in G$ (for instance,
$G$ is finite). Let $\Hom(V,W)$ be the vector space of all linear maps from $V$ to $W$. Using  $\Hom_a(V,W)$ to denote the set of all homogenous linear maps of degree $a$ in $\Hom(V,W)$, we may obtain 
a subspace decomposition of $\Hom(V,W)$ as follows: 
\begin{equation}
\label{ }
\Hom(V,W)=\bigoplus_{a\in G} \Hom_a(V,W).
\end{equation}
This means that $\Hom(V,W)$ is also a $G$-graded vector space. Moreover, if $W=V$, then $\Hom(V,V)$ becomes an associative $G$-graded algebra with respect to the 
usual composition. Clearly, it is also a left-symmetric color algebra for any $\upepsilon\in B_k(G)$. We write 
$\gl_{G,k}(V)$ for this algebra. 
\hbo}\end{exam}

\begin{exam}\label{exam2.9}
{\rm Let $G=\F_2^3, k=\C$, and $\upepsilon\in B_\C(\F_2^3)$ be the second skew-symmetric bicharacter defined in Example \ref{exam2.3}. Let $L$ be a $3$-dimensional $G$-graded vector space over $\C$, spanned by $\{x,y,z\}$, and with
$\deg(x)=(1,1,0),\deg(y)=(1,0,1)$, and $z=(0,1,1)$. In particular, $L_a=0$ if $a\notin\{(1,1,0),(1,0,1),(0,1,1)\}$.
The following nonzero products:
$$xy=z,~~yx=-z$$
define a  left-symmetric color algebra structure on $L$.
\hbo}\end{exam}

\subsection{Lie color algebras} 
A nonassociative $G$-graded algebra $L=\bigoplus_{a\in G}L_a$ with a bracket product $[-,-]$ is called a \textit{Lie color algebra} (or \textit{generalized Lie algebra} in \cite{Sch79}) if there exists an 
$\upepsilon\in B_k(G)$ such that $[x,y]+\upepsilon(|x|,|y|)[y,x]=0$ and
\begin{equation}
\label{ }
\upepsilon(|z|,|x|)[[x,y],z]+\upepsilon(|x|,|y|)[[y,z],x]+\upepsilon(|y|,|z|)[[z,x],y]=0
\end{equation}
for all homogeneous elements $x,y,z\in L$.

Suppose that $A$ denotes a left-symmetric color algebra with respect to $(G,\upepsilon)$. The \textit{$\upepsilon$-commutator}
defined by 
\begin{equation}
\label{ }
[x,y]:=xy-\upepsilon(|x|,|y|)yx
\end{equation}
gives rise to a Lie color algebra structure (with respect to the same $(G,\upepsilon)$) on the underlying vector space of $A$, and denoted as $[A]$. In particular, $[\gl_{G,k}(V)]$ is called  the \textit{general linear Lie color algebra} with respect to $G$ over $k$. 

\begin{exam}\label{exam2.10}
{\rm
Suppose that $G=\F_2^3, k=\C$, and $\upepsilon\in B_\C(\F_2^3)$ be the first skew-symmetric bicharacter defined in Example \ref{exam2.3}. Let $L$ be a $3$-dimensional $G$-graded vector space over $\C$, spanned by $\{x,y,z\}$, and with
$|x|=(1,1,0),|y|=(1,0,1)$, and $|z|=(0,1,1)$. In particular, $L_a=0$ if $a\notin\{(1,1,0),(1,0,1),(0,1,1)\}$.
The following two sets of nonzero bracket products:
$$\{[x,y]=z,[z,x]=y,[y,z]=x\}\textrm{ and }\{[x,y]=z,[z,x]=y,[y,z]=0\}$$
both define a Lie color algebra that is neither a Lie algebra nor a Lie superalgebra. A classification of 3-dimensional
$\F_2^3$-graded complex Lie color algebras with some special conditions has been obtained in \cite{Sil97}.
\hbo}\end{exam}

A $G$-graded vector space $W$ is called a \textit{left module} of a Lie color algebra $L$ (with respect to $(G,\upepsilon)$) if there exists a left action $L\times W\ra W$ such that
\begin{equation}
\label{ }
[x,y]w=x(yw)-\upepsilon(|x|,|y|)y(xw)
\end{equation}
for all homogeneous elements $x,y\in L$, and $w\in W$. Similarly, a \textit{right module} of $L$ is a $G$-graded vector space $W$ with a right action $W\times L\ra W$ such that
\begin{equation}
\label{ }
w[x,y]=(wx)y-\upepsilon(|x|,|y|)(wy)x
\end{equation}
for all homogeneous elements $x,y\in L$, and $w\in W$. 

\subsection{Bimodules}
Let $A$ be a left-symmetric color algebra with  respect to $(G,\upepsilon)$. A $G$-graded vector space $W$ is called an
\textit{$A$-bimodule}  if there exists a left action $A\times W\ra W$ and a right action $W\times A\ra W$ such that
\begin{eqnarray}
(xy)w-x(yw)&=&\upepsilon(|x|,|y|)\Big((yx)w-y(xw)\Big)\label{bm1}\\
(xw)y-x(wy)&=&\upepsilon(|x|,|w|)\Big((wx)y-w(xy)\Big)\label{bm2}
\end{eqnarray}
for all homogeneous elements $x,y\in L$, and $w\in W$.  

\begin{exam}\label{}
{\rm The left and right multiplications of a  left-symmetric color algebra $A$ defines a bimodule structure on $A$ itself, which is called the \textit{natural bimodule} of $A$.
\hbo}\end{exam}

\begin{rem}\label{}
{\rm 
Note that the first condition (\ref{bm1}) above is equivalent to saying that the left action of $A$ on $W$ gives $W$ a left module structure of the Lie color algebra $[A]$. However, the second condition  (\ref{bm2}) doesn't guarantee that the right action of $A$ on $W$ would induce a right module structure of $[A]$ on $W$. 
\hbo}\end{rem}

\begin{defn}\label{}
{\rm Let $A$ be a left-symmetric color algebra.  An $A$-bimodule $W$ is called
\begin{enumerate}
  \item \textit{complete} if the right action of $A$ on $W$ gives a right $[A]$-module structure on $W$;
  \item \textit{right-trivial} if the right action of $A$ on $W$ is trivial, i.e., $wx=0$ for all $x\in A$ and $w\in W$.
\end{enumerate}
}\end{defn}
Clearly, a right-trivial $A$-bimodule must be complete. To construct new bimodules of $A$, let's recall the fact that the tensor product $V\otimes W$ of two $G$-graded vector spaces $V$ and $W$ also can be $G$-graded by
$(V\otimes W)_a=\bigoplus_{b\in G}(V_b\otimes W_{a-b})$, for all $a\in G$.

\begin{prop}\label{prop2.14}
Let $V$ and $W$ be two bimodules of a  left-symmetric color algebra  $A$ and suppose $V$ is complete. Then $V\otimes W$ is an $A$-bimodule with the following left and right actions: 
\begin{eqnarray*}
x(v\otimes w) &:= & (xv-\upepsilon(|x|,|v|)vx)\otimes w+\upepsilon(|x|,|v|)v\otimes (xw) \\
(v\otimes w)x &:= & v\otimes (wx) 
\end{eqnarray*}
for all homogeneous elements $x\in A,v\in V$, and $w\in W$.
\end{prop}

\begin{proof} We assume  that $x\in A_a,y\in A_b,v\in V_c$, and $w\in W_d$ are homogeneous.
To confirm the first condition (\ref{bm1}), we need to show that 
\begin{equation}
\label{ }
[x,y](v\otimes w)=x(y(v\otimes w))-\upepsilon(a,b)y(x(v\otimes w)).
\end{equation}
 In fact, by definition, 
$[x,y](v\otimes w)=\Big([x,y]v-\upepsilon(a+b,c)v[xy]\Big)\otimes w+\upepsilon(a+b,c)v\otimes ([x,y]w)$ and
\begin{eqnarray*}
x(y(v\otimes w))&=&x\Big((yv)\otimes w-\upepsilon(b,c)(vy)\otimes w+\upepsilon(b,c)v\otimes (yw)\Big) \\
&=&\Big(x(yv)-\upepsilon(a,b+c)(yv)x\Big)\otimes w+\upepsilon(a,b+c)(yv)\otimes (xw)\\
&&-\upepsilon(b,c)\Big((x(vy)-\upepsilon(a,b+c)(vy)x)\otimes w+\upepsilon(a,b+c)(vy)\otimes (xw)\Big)\\
&&+\upepsilon(b,c) \Big((xv-\upepsilon(a,c)vx)\otimes yw+\upepsilon(a,c)v\otimes (x(yw))\Big)\\
y(x(v\otimes w))&=&\Big(y(xv)-\upepsilon(b,a+c)(xv)y\Big)\otimes w+\upepsilon(b,a+c)(xv)\otimes (yw)\\
&&-\upepsilon(a,c)\Big((y(vx)-\upepsilon(b,a+c)(vx)y)\otimes w+\upepsilon(b,a+c)(vx)\otimes (yw)\Big)\\
&&+\upepsilon(a,c) \Big((yv-\upepsilon(b,c)vy)\otimes xw+\upepsilon(b,c)v\otimes (y(xw))\Big).
\end{eqnarray*}
Thus $x(y(v\otimes w))-\upepsilon(a,b)y(x(v\otimes w))=([x,y]v)\otimes w-\Big(\upepsilon(a,b+c)(yv)x-\upepsilon(b,c)(xv)y\Big)\otimes w-
\Big(\upepsilon(b,c)x(vy)-\upepsilon(a,b+c)y(vx)\Big)\otimes w-\upepsilon(a+b,c)(v[x,y])\otimes w+\upepsilon(a+b,c)v\otimes ([x,y]w).$
Comparing $[x,y](v\otimes w)$ with $x(y(v\otimes w))-\upepsilon(a,b)y(x(v\otimes w))$,  it suffices to verify that 
\begin{equation}
\label{ }
\upepsilon(a,b+c)\Big((yv)x-y(vx)\Big)=\upepsilon(b,c)\Big((xv)y-x(vy)\Big).
\end{equation}
Note that  $V$ is an $A$-bimodule, it follows that $(xv)y-x(vy)=\upepsilon(a,c)((vx)y-v(xy))$ and $(yv)x-y(vx)=\upepsilon(b,c)((vy)x-v(yx))$. Hence, it suffices to show that 
\begin{equation}
\label{ }
(vx)y-v(xy)=\upepsilon(a,b)((vy)x-v(yx)).
\end{equation}
Actually, this is immediate because $V$ is a complete bimodule, the right action gives a right $[A]$-module structure on $V$, so 
$v(xy-\upepsilon(a,b)yx)=v[x,y]=(vx)y-\upepsilon(a,b)(vy)x$. Therefore, $(vx)y-v(xy)=\upepsilon(a,b)((vy)x-v(yx))$, as desired.

To verify the second condition (\ref{bm2}), we have to show that
\begin{equation}
\label{ }
(x(v\otimes w))y-x((v\otimes w)y)=\upepsilon(a,c+d)\Big(((v\otimes w)x)y-(v\otimes w)(xy)\Big).
\end{equation}
Note that  $(x(v\otimes w))y=\Big((xv-\upepsilon(a,c)vx)\otimes w+\upepsilon(a,c)v\otimes (xw)\Big)y=(xv-\upepsilon(a,c)vx)\otimes (wy)+\upepsilon(a,c)v\otimes ((xw)y)$ and $x((v\otimes w)y)=x(v\otimes wy)=(xv-\upepsilon(a,c)vx)\otimes (wy)+\upepsilon(a,c)v\otimes (x(wy))$. Thus
\begin{eqnarray*}
(x(v\otimes w))y-x((v\otimes w)y)&=&\upepsilon(a,c)v\otimes \Big((xw)y-x(wy)\Big)\\
&=&\upepsilon(a,c)v\otimes \upepsilon(a,d)\Big((wx)y-w(xy)\Big)\\
&=&\upepsilon(a,c+d)\Big(((v\otimes w)x)y-(v\otimes w)(xy)\Big).
\end{eqnarray*}
This completes the proof. 
\end{proof}

\begin{rem}\label{}{\rm 
Note that the $A$-bimodule structure on the tensor product $V\otimes W$ defined above is different from the $A$-bimodule structure in  \cite[Example 2 (iv)]{NBN09}.
\hbo}\end{rem}

The following result gives an $A$-bimodule structure of $\Hom(A,V)$ that actually plays a key role in establishing connections between Lie color cohomology and left-symmetric cohomology in Theorem \ref{maint} below.

\begin{prop}\label{bim1}
Let $V$ be a bimodule of a left-symmetric color algebra $A$ (with respect to $(G,\upepsilon)$) and $\Hom(A,V)$ be the
vector space of all $k$-linear maps from $A$ to $V$, graded by the abelian group $G$.  Together with the trivial right action, the following left action
\begin{equation}
\label{ }
(xf)(z):=xf(z)-\upepsilon(|x|,|f|)f(xz)+\upepsilon(|x|,|f|)f(x)z
\end{equation}
for all homogeneous elements $x,z\in A$ and $f\in \Hom(A,V)$, defines a right-trivial $A$-bimodule structure on $\Hom(A,V)$.
\end{prop}

\begin{proof}
As the right action is trivial, it suffices to verify the condition  (\ref{bm1}). Consider arbitrary homogeneous elements $x\in A_a,y\in A_b,z\in A_d$ and $f\in \Hom_c(A,V)$.
Note that $$([x,y]f)(z)=[x,y]f(z)-\upepsilon(a+b,c)f([x,y]z)+\upepsilon(a+b,c)f([x,y])z$$
and
\begin{eqnarray*}
(x(yf))(z)&=& x\Big((yf)(z)\Big)-\upepsilon(a,b+c)(yf)(xz)+\upepsilon(a,b+c)((yf)(x))z \\
&=&x\Big(yf(z)-\upepsilon(b,c)f(yz)+\upepsilon(b,c)f(y)z\Big)-\\
&&\upepsilon(a,b+c)\Big(yf(xz)-\upepsilon(b,c)f(y(xz))+\upepsilon(b,c)f(y)(xz)\Big)+\\
&&\upepsilon(a,b+c)(yf(x)-\upepsilon(b,c)f(yx)+\upepsilon(b,c)f(y)x)z\\
(y(xf))(z)&=&y\Big(xf(z)-\upepsilon(a,c)f(xz)+\upepsilon(a,c)f(x)z\Big)-\\
&&\upepsilon(b,a+c)\Big(xf(yz)-\upepsilon(a,c)f(x(yz))+\upepsilon(a,c)f(x)(yz)\Big)+\\
&&\upepsilon(b,a+c)(xf(y)-\upepsilon(a,c)f(xy)+\upepsilon(a,c)f(x)y)z.
\end{eqnarray*}
Using the two facts that $x(f(y)z)-(xf(y))z=\upepsilon(a,b+c)[f(y)(xz)-(f(y)x)z]$ and 
$y(f(x)z)-(yf(x))z=\upepsilon(b,a+c)[f(x)(yz)-(f(x)y)z]$, a direct computation shows that
\begin{eqnarray*}
\Big(x(yf)-\upepsilon(a,b)y(xf)\Big)(z)& = & (x(yf))(z)-\upepsilon(a,b)(y(xf))(z) \\
 & = & [x,y]f(z)-\upepsilon(a+b,c)f([x,y]z)+\upepsilon(a+b,c)f([x,y])z\\
 &=& ([x,y]f)(z).
\end{eqnarray*}
Hence, $x(yf)-\upepsilon(a,b)y(xf)=[x,y]f$ and the the condition  (\ref{bm1}) is confirmed.
\end{proof}

\section{Cohomology} \label{sec3}
\setcounter{equation}{0}
\renewcommand{\theequation}
{3.\arabic{equation}}
\setcounter{theorem}{0}
\renewcommand{\thetheorem}
{3.\arabic{theorem}}

\subsection{Coboundary operators}
Let $\upepsilon\in B_k(G)$ and $A$ be a left-symmetric color algebra with respect to $\upepsilon$.
Suppose $V$ is a finite-dimensional $A$-bimodule. To construct the cohomology space of $A$ with coefficients in $V$, we define 
\begin{equation}
\label{ }
C^n(A,V):=\Hom((\land^{n-1}_\upepsilon A)\otimes A,V)=\bigoplus_{c\in G}\Hom_c\left((\land^{n-1}_\upepsilon A)\otimes A,V\right)
\end{equation}
as the space of \textit{$n$-th cochains}, which can be $G$-graded.  Here $\land^{n-1}_\upepsilon A$ denotes the $(n-1)$-th $\upepsilon$-exterior power of $A$, which is also a $G$-graded vector space, and $\Hom_c\left((\land^{n-1}_\upepsilon A)\otimes A,V\right)$
denotes the vector space of all homogeneous linear maps of degree $c$ from $(\land^{n-1}_\upepsilon A)\otimes A$ to $V$; see \cite[Appendix A]{SZ98} for more details about $\upepsilon$-exterior algebras. 

In particular, $C^1(A,V)=\Hom(A,V)$ and we define the \textit{$n$-th coboundary operator}
$d_n:C^n(A,V)\ra C^{n+1}(A,V)$ by
\begin{eqnarray*}
\label{ }
(d_nf)(x_1,\dots,x_{n+1})&:=&\sum_{i=1}^n(-1)^{i+1}\upepsilon\left(|f|+\sum_{j=1}^{i-1}|x_j|,|x_i|\right)x_i f(x_1,\dots,\ha{x_i},\dots,x_{n+1})\\
&&\hspace{-0.4cm}+\sum_{i=1}^n(-1)^{i+1}\upepsilon\left(|x_i|, \sum_{j=i+1}^n|x_j|\right)f(x_1,\dots,\ha{x_i},\dots,x_n,x_i)x_{n+1}\\
&&\hspace{-0.4cm}-\sum_{i=1}^n(-1)^{i+1}\upepsilon\left(|x_i|, \sum_{j=i+1}^n|x_j|\right)f(x_1,\dots,\ha{x_i},\dots,x_n,x_ix_{n+1})\\
&&\hspace{-0.4cm}+\sum_{1\leqslant j<i\leqslant n}(-1)^{i+1}\upepsilon\left(\sum_{s=j+1}^{i-1}|x_s|,|x_i|\right)f(x_1,\dots,[x_j,x_i],x_{j+1},\dots,\ha{x_i},\dots,x_{n+1})
\end{eqnarray*}
for all homogeneous elements $f\in C^n(A,V)$ and $x_1,\dots,x_{n+1}\in A$.

\begin{rem}\label{}
{\rm 
Our cochain space $C^n(A,V)$ and the coboundary operators $d_n$ both are different from  \cite[Section 3.1]{NBN09}.
\hbo}\end{rem}

\begin{exam}\label{exam3.2}
{\rm
The first coboundary operator $d_1:C^1(A,V)\ra C^2(A,V)$ is given by
\begin{equation}
\label{ }
(d_1f)(x_1,x_{2})=\upepsilon\left(c,a\right)x_1f(x_2)+f(x_1)x_2-f(x_1x_2)
\end{equation}
where $|x_1|=a,|x_2|=b$, and $|f|=c$. The second coboundary operator $d_2:C^2(A,V)\ra C^3(A,V)$ is defined by
\begin{eqnarray*}
(d_2f)(x_1,x_{2},x_3)&=&\upepsilon(d,a)x_1f(x_2,x_3)-\upepsilon(a+d,b)x_2f(x_1,x_3)+\upepsilon(a,b)f(x_2,x_1)x_3\\
&&-f(x_1,x_2)x_3-\upepsilon(a,b)f(x_2,x_1x_3)+f(x_1,x_2x_3)-f([x_1,x_2],x_3)
\end{eqnarray*}
where $|x_1|=a,|x_2|=b, |x_3|=c$, and $|f|=d$. Substituting $f$ with $d_1f$ in the previous equation obtains
\begin{eqnarray*}
&&((d_2\circ d_1)(f))(x_1,x_{2},x_3)=(d_2(d_1f))(x_1,x_{2},x_3)\\
&=&\upepsilon(d,a)x_1(d_1f)(x_2,x_3)-\upepsilon(a+d,b)x_2(d_1f)(x_1,x_3)+\upepsilon(a,b)(d_1f)(x_2,x_1)x_3\\
&&-(d_1f)(x_1,x_2)x_3-\upepsilon(a,b)(d_1f)(x_2,x_1x_3)+(d_1f)(x_1,x_2x_3)-(d_1f)([x_1,x_2],x_3)
\end{eqnarray*}
which is equal to zero. In fact,
\begin{eqnarray*}
\upepsilon(d,a)x_1(d_1f)(x_2,x_3)&=&\upepsilon(d,a)x_1[\upepsilon\left(d,b\right)x_2f(x_3)+f(x_2)x_3-f(x_2x_3)]\\
-\upepsilon(a+d,b)x_2(d_1f)(x_1,x_3)&=&-\upepsilon(a+d,b)x_2[\upepsilon\left(d,a\right)x_1f(x_3)+f(x_1)x_3-f(x_1x_3)]\\
\upepsilon(a,b)(d_1f)(x_2,x_1)x_3&=&\upepsilon(a,b)[\upepsilon\left(d,b\right)x_2f(x_1)+f(x_2)x_1-f(x_2x_1)]x_3\\
-(d_1f)(x_1,x_2)x_3&=&-[\upepsilon\left(d,a\right)x_1f(x_2)+f(x_1)x_2-f(x_1x_2)]x_3\\
-\upepsilon(a,b)(d_1f)(x_2,x_1x_3)&=&-\upepsilon(a,b)[\upepsilon\left(d,b\right)x_2f(x_1x_3)+f(x_2)(x_1x_3)-f(x_2(x_1x_3))]\\
(d_1f)(x_1,x_2x_3)&=&\upepsilon\left(d,a\right)x_1f(x_2x_3)+f(x_1)(x_2x_3)-f(x_1(x_2x_3))\\
-f([x_1,x_2],x_3)&=&-[\upepsilon\left(d,a+b\right)[x_1,x_2]f(x_3)+f([x_1,x_2])x_3-f([x_1,x_2]x_3)]
\end{eqnarray*}
Using (\ref{lsa}), (\ref{bm1}), and (\ref{bm2}) to the sum of all right-hand sides of these equations,
we may see that the sum could be simplified to zero. This means that $d_2\circ d_1=0.$
\hbo}\end{exam}

\begin{rem}{\rm
We make the convention that $C^0(A,V):=\{v\in V\mid (xy)v=x(yv),\textrm{ for all }x,y\in A\}$ and 
extend the notion $C^n(A,V)$ for all $n\in\N$. It is clear to see that  $C^0(A,V)$ is an $A$-bimodule and so it is a
submodule of $V$.  We define the $0$-th coboundary operator 
$d_0:C^0(A,V)\ra C^1(A,V)$ by
\begin{equation}
\label{ }
d_0(v)(x):=vx-\upepsilon(|v|,|x|)xv
\end{equation}
for all homogeneous elements $v\in C^0(A,V)$ and $x\in A$. Moreover, using (\ref{bm2}) and the definition of $C^0(A,V)$, we 
see that
\begin{eqnarray*}
(d_1\circ d_0)(v)(x_1,x_{2}) & = & d_1(d_0(v))(x_1,x_{2}) \\
 & = & \upepsilon\left(c,a\right)x_1\cdot d_0(v)(x_2)+d_0(v)(x_1)\cdot x_2-d_0(v)(x_1x_2) \\
 &=&\upepsilon\left(c,a\right)x_1(vx_2-\upepsilon(c,b)x_2v)+(vx_1-\upepsilon(c,a)x_1v)x_2-(v(x_1x_2)\\
 &&-\upepsilon(c,a+b)(x_1x_2)v)\\
 &=&0
\end{eqnarray*}
where $|x_1|=a,|x_2|=b,$ and $|v|=c$. Hence, we obtain a chain of $G$-graded vector spaces:
\begin{equation}
\label{cochain}
\xymatrix{
C^0(A,V) \ar[r]^-{d_0} & C^1(A,V) \ar[r]^-{d_1}  & C^2(A,V) \ar[r]^-{d_2} & C^3(A,V) \ar[r]^-{d_3} & \cdots 
}
\end{equation}
with $d_1\circ d_0=d_2\circ d_1=0.$
\hbo}\end{rem}

More generally,

\begin{prop}\label{cob}
In (\ref{cochain}), we have $d_{n+1}\circ d_n=0$  for all $n\in\N$.
\end{prop}

\begin{rem}{\rm
This statement could be verified by a direct but tedious computation via repeatedly applying (\ref{bm1}), (\ref{bm2}), and (\ref{lsa}).
However, we will see that Lemma \ref{lem2} below, together with the fact in Lie color cohomology that the composition of any two adjacent coboundary operators is zero, can provide a shorter proof for Proposition \ref{cob}.
\hbo}\end{rem}

\begin{rem}{\rm
Recall that $C^1(A,V)=\Hom(A,V)$. We actually are able to generalize Proposition \ref{bim1} to the case of all 
$C^n(A,V)$. More specifically, the following left action
\begin{eqnarray*}
(xf)(x_1,\dots,x_{n+1})&:=&xf(x_1,\dots,x_{n+1})-\upepsilon\left(|x|,|f|+\sum_{i=1}^n|x_i|\right)f(x_1,\dots,x_n,xx_{n+1})\\
&&-\sum_{j=1}^n \upepsilon\left(|x|,|f|+\sum_{i=1}^{j-1}|z_i|\right)f(x_1,\cdots,x_{j-1},[x,x_j], x_{j+1},\cdots,x_n,
x_{n+1})\\
&&+\upepsilon\left(|x|,|f|+\sum_{i=1}^n|x_i|\right)f(x_1,\dots,x_{n+1})x
\end{eqnarray*}
together with the trivial right action, defines a right-trivial $A$-module structure on the space $C^{n+1}(A,V)$ for all $n\in\N$.
\hbo}\end{rem}

\subsection{Cohomology}
Suppose $A$ denotes a left-symmetric color algebra with respect to $\upepsilon$ and $V$ is a finite-dimensional 
$A$-bimodule. For all $n\in\N$, we call
\begin{equation}
\label{ }
Z^n(A,V):=\{f\in C^n(A,V)\mid d_n(f)=0\}
\end{equation}
the space of  \textit{$n$-th cocycles} of $A$ with coefficients in $V$. For $n\in\N^+$, we use $B^n(A,V)$ to denote the image of $d_{n-1}$ in $C^{n}(A,V)$, 
which is a $G$-graded subspace of $C^{n}(A,V)$ and called the space of \textit{$n$-th coboundaries} of $A$ with coefficients in $V$. 
By Proposition \ref{cob}, we see that $B^{n}(A,V)\subseteq Z^{n}(A,V)$. 
Note that $Z^n(A,V)$ is also a $G$-graded subspace of $C^n(A,V)$. The quotient space
\begin{equation}
\label{sim}
H^n(A,V):=\frac{Z^{n}(A,V)}{B^{n}(A,V)}
\end{equation}
is called the \textit{$n$-th cohomology space} of $A$ with coefficients in $V$, for all $n\geqslant 1$. 

Note that the coboundary operators $d_n$ are degree-preserving, thus $B^{n}(A,V)$ is a $G$-graded subspace of $Z^{n}(A,V)$, i.e.,
$B_c^{n}(A,V)=Z_c^{n}(A,V)\cap B^{n}(A,V)$, for all $c\in G$. Therefore, 
$$H^n(A,V)=\bigoplus_{c\in G}H^n_c(A,V)$$
is also a $G$-graded vector space, where
\begin{equation}
\label{coh}
H^n_c(A,V)=\frac{Z^n_c(A,V)+B^{n}(A,V)}{B^{n}(A,V)}\cong \frac{Z^n_c(A,V)}{Z^n_c(A,V)\cap B^{n}(A,V)}=\frac{Z^n_c(A,V)}{B^{n}_c(A,V)}.
\end{equation}
Clearly, for each $c\in G$, $H^n_c(A,V)=0$ if and only if $Z^n_c(A,V)=B^n_c(A,V)$.

\begin{rem}\label{rem3.7}
{\rm
We define the \textit{$0$-th cohomology space} $H^0(A,V)$ to be
\begin{eqnarray*}
H^0(A,V)&:=&Z^{0}(A,V)=\{v\in C^0(A,V)\mid d_0(v)=0\}\\
&=&\{v\in V\mid vx=\upepsilon(|v|,|x|)xv\textrm{ and }(xy)v=x(yv),\textrm{ for all }x,y\in A\}\\
&=&\{v\in V\mid vx=\upepsilon(|v|,|x|)xv,\textrm{ for all }x\in A\}
\end{eqnarray*}
where the condition ``$vx=\upepsilon(|v|,|x|)xv$ for all $x\in A$'' is sufficient for the condition 
``$(xy)v=x(yv),\textrm{ for all }x,y\in A$''. Moreover, assume that $V=A$ denotes the natural bimodule. Then
\begin{eqnarray*}
H^0(A,A)&=&\{v\in A\mid vu-\upepsilon(|v|,|u|)uv=0,\textrm{ for all }u\in A\}\\
&=&\{v\in A\mid [u,v]=0,\textrm{ for all }u\in [A]\}\\
&=&A^{[A]}
\end{eqnarray*}
is the subspace of invariants of $[A]$ on itself. Thus, the $0$-th cohomology space $H^0(A,V)$ is also called the \textit{space of invariants} of $A$ on $V$. 

Let's look back at Example \ref{exam2.9}.  The space $H^0(A,A)$ can be computed by 
$$\{v\in A\mid vx-\upepsilon(|v|,|x|)xv=0\}\cap \{v\in A\mid vy-\upepsilon(|v|,|y|)yv=0\}\cap \{v\in A\mid vz-\upepsilon(|v|,|z|)zv=0\},$$
in which the three subspaces are spanned by $\{x,z\}$, $\{y,z\}$, and $\{x,y,z\}$, respectively.  Hence, $H^0(A,A)$ is a one-dimensional $G$-graded vector space spanned by $\{z\}$, i.e., 
$$H^0(A,A)=H^0_{(0,1,1)}(A,A)=\C\cdot z$$ and $H^0_{a}(A,A)=0$ for all $a\in G$ but $a\neq (0,1,1)$.
\hbo}\end{rem}

Furthermore, together with \cite[Theorem 3.2]{Fel01}, Remark \ref{rem3.7} implies the following result about $H^0(A,A)$ for some special left-symmetric color algebras $A$.

\begin{coro}
If the left multiplication $\ell_x: A\ra A$ defined by $y\mapsto xy$ is nilpotent for each homogeneous $x\in A$, then
$H^0(A,A)\neq 0.$ 
\end{coro}

\begin{exam}\label{coh1}
{\rm
The coboundary space $B^1(A,V)$ consists of all homogeneous linear maps $f$ for which there exists a homogeneous element $v\in V$
such that $f(x)=vx-\upepsilon(|v|,|x|)xv$ for all homogeneous $x\in A$. In particular, if $V=A$ denotes the natural bimodule, 
each element of $B^1(A,A)$ can be regarded as an inner derivation of the Lie color algebra $[A]$. 

The cocycle space $Z^1(A,V)$ consists of all homogeneous linear maps $f$ satisfying 
\begin{equation}
\label{ }
f(x_1x_2)=f(x_1)x_2+\upepsilon\left(|f|,|x_1|\right)x_1f(x_2),
\end{equation}
which is exactly equal to the space of $\upepsilon$-derivations of $A$ to $V$ in terms of \cite[Examples 2.4]{Sch83}.

Thus it is reasonable to call each element of $B^1(A,V)$ an \textit{inner $\upepsilon$-derivation} of $A$ to $V$, and call the first cohomology space $H^1(A,V)$ the \textit{space of outer $\upepsilon$-derivations} of $A$ to $V$.
\hbo}\end{exam}

\subsection{Connection with Lie color cohomology}

Let $L$ be a Lie color algebra with respect to $\upepsilon\in B_k(G)$ and $W$ be a left $L$-module. 
We still write $C^n(L,W)$ for the space of $n$-cochains of $L$ in $W$, which is defined as the $G$-graded vector space of all $\upepsilon$-alternating $n$-multilinear linear mappings on $L$ to $W$. Serving as a fundamental object in defining Lie color cocycles, Lie color coboundaries, and Lie color cohomology, the space $C^n(L,W)$ can be identified with the $G$-graded space $$\Hom(\land^n_\upepsilon L,W)=\bigoplus_{c\in G}\Hom_c(\land^n_\upepsilon L,W).$$ In particular, for $n\in\N^+$,  the $n$-th coboundary operator $\delta_n: C^n(L,W)\ra C^{n+1}(L,W)$ is given by
\begin{eqnarray*}
\delta_n(f)(x_1,\dots,x_{n+1})&:=&\sum_{i=1}^{n+1}(-1)^{i+1}\upepsilon \left(|f|+\sum_{j=1}^{i-1}|x_j|,|x_i|\right) x_i f(x_1,\dots,\ha{x_{i}},\dots,x_{n+1})\\
&&\hspace{-3cm}+\sum_{1\leqslant j<i\leqslant n+1} (-1)^{i+1} \upepsilon \left(\sum_{s=j+1}^{i-1}|x_s|,|x_i|\right) f(x_1,\dots,x_{j-1},[x_j,x_i],x_{j+1},\dots,\ha{x_{i}},\dots,x_{n+1})
\end{eqnarray*}
for all homogenous elements $x_1,\dots,x_{n+1}\in L$ and $f\in C^n(L,W)$. See for exmaple \cite[Section 2]{SZ98} for more information.

\begin{lem} \label{lem1}
Let $V$ be a bimodule of a left-symmetric color algebra $A$ with respect to $\upepsilon$. Then there exists an isomorphism 
of  $G$-graded vector spaces: 
$$C^{n+1}(A,V)\cong C^n([A],C^1(A,V))$$
for all $n\in\N$. 
\end{lem}

\begin{proof} 
By the graded tensor-hom adjunction (see for example \cite[Proposition 6.15]{Lop02}), there exists a natural isomorphism of  $G$-graded vector spaces:
$$C^{n+1}(A,V)=\Hom((\land^n_\upepsilon A)\otimes A,V)\cong\Hom(\land^n_\upepsilon A,\Hom(A,V))=
\Hom(\land^n_\upepsilon [A],\Hom(A,V)).$$
It follows from Proposition \ref{bim1} that  $\Hom(A,V)=C^1(A,V)$ is a left $[A]$-module.
Thus $$\Hom(\land^n_\upepsilon [A],\Hom(A,V))=\Hom(\land^n_\upepsilon [A],C^1(A,V))=C^n([A],C^1(A,V)).$$
Therefore, $C^{n+1}(A,V)$ and $C^n([A],C^1(A,V))$ are isomorphic as $G$-graded vector spaces. 
\end{proof}

\begin{rem}{\rm
The natural isomorphism $\varphi: C^{n+1}(A,V)\ra C^n([A],C^1(A,V))$ can be described  explicitly by
\begin{equation}
\label{ni}
\Big(\varphi(f)(x_1,\dots,x_n)\Big)(x):=f(x_1,\dots,x_n,x)
\end{equation}
for all $x_1,\dots,x_n,x\in A$ and $f\in C^{n+1}(A,V)$.
\hbo}\end{rem}

\begin{lem}\label{lem2}
Let $\varphi: C^{n+1}(A,V)\ra C^n([A],C^1(A,V))$ be the natural isomorphism (\ref{ni}) and
$\delta_n: C^n([A],C^1(A,V))\ra C^{n+1}([A],C^1(A,V))$ be the coboundary operator. Then the following diagram commutes:
$$\xymatrix{
C^{n+1}(A,V) \ar[rr]^-{\varphi}\ar[d]^-{d_{n+1}} && C^n([A],C^1(A,V)) \ar[d]_-{\delta_{n}}  \\
C^{n+2}(A,V) \ar[rr]^-{\varphi} && C^{n+1}([A],C^1(A,V)). 
}$$
\end{lem}

\begin{proof}
Let $f\in C^{n+1}(A,V), x_1,\dots,x_{n+1}\in [A]$, and $x\in A$ be arbitrary homogeneous elements.  We need to show that
$(\delta_n\circ \varphi)(f)=(\varphi\circ d_{n+1})(f)$. In fact, the map $(\delta_n\circ \varphi)(f)$ is given by
\begin{eqnarray*}
&&\Big((\delta_n\circ \varphi)(f)(x_1,\dots,x_{n+1})\Big)(x)= \Big((\delta_n(\varphi(f))(x_1,\dots,x_{n+1})\Big)(x)  \\
 &=&\sum_{i=1}^{n+1}(-1)^{i+1}\upepsilon \left(|f|+\sum_{j=1}^{i-1}|x_j|,|x_i|\right) \Big(x_i (\varphi(f)(x_1,\dots,\ha{x_{i}},\dots,x_{n+1}))\Big)(x)+\\
&&\hspace{-1cm}\sum_{1\leqslant j<i\leqslant n+1} (-1)^{i+1} \upepsilon \left(\sum_{s=j+1}^{i-1}|x_s|,|x_i|\right) \Big(\varphi(f)(x_1,\dots,x_{j-1},[x_j,x_i],x_{j+1},\dots,\ha{x_{i}},\dots,x_{n+1})\Big)(x).
\end{eqnarray*}
Note that $\varphi(f)(x_1,\dots,\ha{x_{i}},\dots,x_{n+1})\in C^1(A,V)$. By Proposition \ref{bim1}, we have
\begin{eqnarray*}
&&\Big(x_i (\varphi(f)(x_1,\dots,\ha{x_{i}},\dots,x_{n+1}))\Big)(x)= x_i \Big(\varphi(f)(x_1,\dots,\ha{x_{i}},\dots,x_{n+1})(x)\Big)\\
 &&\hspace{2cm} -\upepsilon\left(|x_i|,|f|+\sum_{j=1}^{i-1}|x_j|+\sum_{j=i+1}^{n+1}|x_j|\right)\varphi(f)(x_1,\dots,\ha{x_{i}},\dots,x_{n+1})(x_ix)\\
 &&\hspace{2cm}+\upepsilon\left(|x_i|,|f|+\sum_{j=1}^{i-1}|x_j|+\sum_{j=i+1}^{n+1}|x_j|\right)\Big(\varphi(f)(x_1,\dots,\ha{x_{i}},\dots,x_{n+1})(x_i)\Big)x\\
 &&\hspace{0.5cm}= x_i f(x_1,\dots,\ha{x_{i}},\dots,x_{n+1},x)-\upepsilon\left(|x_i|,|f|+\sum_{j=1}^{i-1}|x_j|+\sum_{j=i+1}^{n+1}|x_j|\right)f(x_1,\dots,\ha{x_{i}},\dots,x_{n+1},x_ix)\\
 &&\hspace{0.9cm}+\upepsilon\left(|x_i|,|f|+\sum_{j=1}^{i-1}|x_j|+\sum_{j=i+1}^{n+1}|x_j|\right)f(x_1,\dots,\ha{x_{i}},\dots,x_{n+1},x_i)x
\end{eqnarray*}
for all $i=1,2,\dots,n+1$. Hence,
\begin{eqnarray*}
&&\Big((\delta_n\circ \varphi)(f)(x_1,\dots,x_{n+1})\Big)(x)\\
 &=&\sum_{i=1}^{n+1}(-1)^{i+1}\upepsilon \left(|f|+\sum_{j=1}^{i-1}|x_j|,|x_i|\right) x_i f(x_1,\dots,\ha{x_{i}},\dots,x_{n+1},x)\\
 &&+\sum_{i=1}^{n+1}(-1)^{i+1} \upepsilon\left(|x_i|,\sum_{j=i+1}^{n+1}|x_j|\right)f(x_1,\dots,\ha{x_{i}},\dots,x_{n+1},x_i)x\\
 &&-\sum_{i=1}^{n+1}(-1)^{i+1}  \upepsilon\left(|x_i|,\sum_{j=i+1}^{n+1}|x_j|\right)f(x_1,\dots,\ha{x_{i}},\dots,x_{n+1},x_ix)\\
&&+\sum_{1\leqslant j<i\leqslant n+1} (-1)^{i+1} \upepsilon \left(\sum_{s=j+1}^{i-1}|x_s|,|x_i|\right)f(x_1,\dots,x_{j-1},[x_j,x_i],x_{j+1},\dots,\ha{x_{i}},\dots,x_{n+1},x)\\
&=&(d_{n+1}f)(x_1,\dots,x_{n+1},x)=\varphi\Big(d_{n+1}(f)(x_1,\dots,x_{n+1})\Big)(x).
\end{eqnarray*}
The last equation holds from the definition of coboundary operators $d_n$.  This shows that 
$(\delta_n\circ \varphi)(f)=(\varphi\circ d_{n+1})(f)$ and therefore, the diagram is commutative.
\end{proof}

\begin{rem}\label{inv}
{\rm
We also have the following commutative diagram for the inverse of $\varphi$:
$$\xymatrix{
C^n([A],C^1(A,V)) \ar[d]_-{\delta_{n}} \ar[rr]^-{\varphi^{-1}}&& C^{n+1}(A,V) \ar[d]_-{d_{n+1}}  \\
C^{n+1}([A],C^1(A,V))\ar[rr]^-{\varphi^{-1}} && C^{n+2}(A,V). 
}$$
In fact, for each $g\in C^n([A],C^1(A,V))$, there exists $f\in C^{n+1}(A,V)$ such that $g=\varphi(f)$ because $\varphi$ is an isomorphism. Thus $\varphi^{-1}(\delta_n(g))=\varphi^{-1}(\delta_n(\varphi(f)))=
\varphi^{-1}(\varphi(d_{n+1}(f)))=d_{n+1}(f)=d_{n+1}(\varphi^{-1}(g))$. Namely, $\varphi^{-1}\circ\delta_n=d_{n+1}\circ\varphi^{-1}$.
\hbo}\end{rem}

Lemma \ref{lem2} can be used to give a short proof for Proposition \ref{cob}.

\begin{proof}[Proof of Proposition \ref{cob}]
We have seen from Example  \ref{exam3.2} that $d_2\circ d_1=0$, so we may suppose that $n\geqslant 2$.
Since $\delta_{n-1}\circ \varphi=\varphi\circ d_{n}$ and $\varphi$ is an isomorphism, it follows that
$d_n=\varphi^{-1}\circ\delta_{n-1}\circ \varphi$. Note that  $\delta_{n}\circ \varphi=\varphi\circ d_{n+1}$, thus
$(\varphi\circ d_{n+1})\circ d_n  =  (\delta_{n}\circ \varphi)\circ d_n=(\delta_{n}\circ \varphi)\circ  (\varphi^{-1}\circ\delta_{n-1}\circ \varphi)=\delta_{n}\circ\delta_{n-1}\circ \varphi=0\circ \varphi=0$. However, as $\varphi$ is isomorphic, we see that 
$d_{n+1}\circ d_n  =0$. Here we have used the fact that $\delta_{n}\circ\delta_{n-1}=0$, because $\delta_{n}$ and $\delta_{n-1}$ are coboundary operators in the Lie color cohomology and their composition $\delta_{n}\circ\delta_{n-1}$ must be zero; see \cite[Section 2]{SZ98}.
\end{proof}

Now we are ready to prove the main result. 

\begin{proof}[Proof of Theorem \ref{maint}]
It suffices to show that $H_c^{n+1}(A,V)\cong H_c^n([A], C^1(A,V))$ for each $c\in G$. Let's take any element $f\in Z_c^{n+1}(A,V)$ of degree $c\in G$. Note that $d_{n+1}(f)=0$ and so it follows from Lemma \ref{lem2} that
$\delta_n(\varphi(f))=\varphi(d_{n+1}(f))=\varphi(0)=0$. Thus $\varphi(f)\in Z_c^n([A], C^1(A,V))$.
As $\varphi$ is an isomorphism, restricting on $Z_c^{n+1}(A,V)$ gives an embedding of vector spaces to 
$Z_c^n([A], C^1(A,V))$. Similarly, by using Remark \ref{inv}, the inverse of $\varphi$ gives rise to an embedding of vector spaces from
$Z_c^n([A], C^1(A,V))$ to $Z_c^{n+1}(A,V)$. Thus 
$$\varphi(Z_c^{n+1}(A,V))=Z_c^n([A], C^1(A,V)).$$ This fact, together with the standard map $\pi: Z_c^n([A], C^1(A,V))\ra H_c^n([A], C^1(A,V))$ obtained from (\ref{coh}), establishes 
a surjective linear map $\theta: Z_c^{n+1}(A,V)\ra H_c^n([A], C^1(A,V))$ defined by $\theta:=\pi\circ \varphi$.

To show that $\ker(\theta)=B_c^{n+1}(A,V)$, we may assume that $f\in B_c^{n+1}(A,V)$. Then there exists
$g\in Z_c^{n}(A,V)$ such that $f=d_n(g)$. By Lemma \ref{lem2}, we see that
$$\varphi(f)=\varphi(d_n(g))=\delta_{n-1}(\varphi(g))\in B_c^n([A], C^1(A,V))$$
and so $\theta (f)=\pi(\varphi(f))=0$, which means that $B_c^{n+1}(A,V)\subseteq \ker(\theta)$. 
Note that $\varphi$ is injective, thus $B_c^{n+1}(A,V)$ can be embedded into $B_c^n([A], C^1(A,V))$ via $\varphi$. 
Conversely, using $\varphi^{-1}$ and Remark \ref{inv} obtains that $B_c^n([A], C^1(A,V))$ can be viewed as
a subspace of $B_c^{n+1}(A,V)$. Thus 
$\varphi(B_c^{n+1}(A,V))=B_c^n([A], C^1(A,V)).$
Now to prove that $\ker(\theta)=B_c^{n+1}(A,V)$, we only need to show that their dimensions are equal to each other. 
In fact, since $\varphi$ is an isomorphism, it follows that $\dim(\ker(\theta))=\dim(\ker(\pi))=\dim(B_c^n([A], C^1(A,V)))=\dim(B_c^{n+1}(A,V))$. Therefore, it follows from (\ref{sim}) that
$$H_c^{n+1}(A,V)\cong\frac{Z_c^{n+1}(A,V)}{B_c^{n+1}(A,V)}\cong H_c^n([A], C^1(A,V))$$
are isomorphic as vector spaces. 
\end{proof}

\begin{rem}{\rm
The first cohomology $H^{1}(A,V)$ can be computed directly as in Example \ref{coh1}.
\hbo}\end{rem}

\begin{rem}{\rm
Given a Lie color algebra $L$ with respect to $(G,\upepsilon)$, we use $\A_L$ to denote the set of all left-symmetric color algebras $A$  with respect to the same $(G,\upepsilon)$ such that $L=[A]$. Suppose $\A_L\neq\emptyset$ and $V$ denotes a bimodule of an algebra $A\in\A_L$. Theorem \ref{maint} allows us to compute the $(n+1)$-th cohomology of $A$ in $V$ by computing 
the $n$-th cohomology of $L$ in $\Hom(A,V)$. 

Moreover, the close connections between Lie color algebras and left-symmetric color algebras lead to the following natural questions:
 (1) When is $\A_L$ not empty? (2) If $\A_L\neq\emptyset$, how can we classify those  left-symmetric color algebras in $\A_L$ up to isomorphism? Similar questions for left-symmetric algebras and left-symmetric algebras have been studied extensively; see for example, \cite{Bau99} and \cite{DZ23}.
\hbo}\end{rem}

\section{Low-dimensional Algebras} \label{sec4}
\setcounter{equation}{0}
\renewcommand{\theequation}
{4.\arabic{equation}}
\setcounter{theorem}{0}
\renewcommand{\thetheorem}
{4.\arabic{theorem}}

\noindent In the study of any nonassociative algebras, computations of low-dimensional cases always play an important role in understanding structures and representations of these algebras; see for example, \cite{BM01, BMH03, CLZ14, CZ17}, and \cite{CZZZ18}. This section is devoted to a discussion on the varieties of two-dimensional and three-dimensional left-symmetric color algebras.

\subsection{$2$-dimensional left-symmetric color algebras}

Suppose that $A=\oplus_{a\in G}A_a$ denotes a nontrivial 2-dimensional left-symmetric color algebra over a field $k$, graded by a finite abelian group $G$ and spanned by $\{x,y\}$. Clearly, $\dim_k(A)=\sum_{a\in G}\dim_k(A_a)$. We may assume that $\upepsilon\in B_k(G)$ and $|G|\geqslant 2$ as if $|G|=1$, then $A$ becomes a left-symmetric algebra, and there exists a classification of 2-dimensional left-symmetric algebras; see \cite[Section 5]{ZB12}. 

The fact that $A_{a}\cdot A_{b}\subseteq A_{a+b}$ is greatly helpful in the analyzing the variety of 2-dimensional left-symmetric color algebras. Note that determining the structure of $A$ is equivalent to finding the values of four products:
$$x^2,~~xy,~~yx,~~y^2.$$

\textbf{Case 1}.  Let's discuss the first case where $x\in A_a, y\in A_b$ for some $a,b\in G\setminus\{0\}$ and $a\neq b$.
Since $\dim(A)=2$ and $\dim(A_a)+\dim(A_b)\geqslant 2$, it follows that $A_{c}=0$ for all $c\notin\{a,b\}$.
As $A$ is nonzero, at least one of $\{x^2,xy,yx,y^2\}$ is not zero, which means that at least one of $\{2a,a+b,2b\}$ 
is equal to either $a$ or $b$. However, note that $|xy|=|yx|=|x|+|y|=a+b$, and both $a$ and $b$ are not zero, thus 
$xy=yx=0$. Therefore, at least one of $\{2a, 2b\}$ belongs to $\{a,b\}$.

\textsc{Subcase 1}. Suppose $2a\in\{a,b\}$ and $2b\notin\{a,b\}$. Note that $a\neq 0$, thus $2a\neq a$ and so $2a=b$. Thus
$A=A_a\oplus A_b$ with $A_a\cdot A_a\subseteq A_b$ and we may assume that $x^2=c_1y$ for some $c_1\in k$.
The assumption that $2b\notin\{a,b\}$ means that $y^2=0$. Combining $xy=yx=0, y^2=0$, and $x^2=c_1y$
into (\ref{lsa}) implies that $A$ must be trivial, i.e., the product of any two elements in $A$ is zero.

\textsc{Subcase 2}. Suppose $2b\in\{a,b\}$ and $2a\notin\{a,b\}$. Symmetrically, switching the roles of $x$ and $y$, we also will obtain a trivial left-symmetrical color algebra $A$.

\textsc{Subcase 3}. Suppose $2a,2b\in\{a,b\}$. Since $2a\neq a$ and $2b\neq b$, we see that
$2a=b$ and $2b=a$. Assume that $x^2=c_1 y$ and $y^2=c_2 x$ some $c_1,c_2\in k$. Setting $z=y$ in 
(\ref{lsa}) we have
$$(xy)y-x(yy)=\upepsilon(|x|,|y|)\Big((yx)y-y(xy)\Big)=0.$$
Thus $0=x(y^2)=x(c_2 x)=c_2 x^2=c_1c_2 y$, and $c_1c_2=0$. Without loss of generality, we assume that $c_2=0$
and obtain a left-symmetric color algebra $A$ defined by the following nonzero product:
\begin{equation}
\label{ }
x^2=c y
\end{equation}
where $c\in k$, and $A=k\cdot x\oplus k\cdot y$ is graded by $G$ for which $|x|+|x|=|y|$.

\textbf{Case 2}.  Consider the second case where $x,y\in A_a$ for some $a\in G\setminus\{0\}$. In this case, $A=A_a$ is spanned by $\{x,y\}$. Note that $|x^2|=2a$ and if $x^2\neq 0$, then $x^2$ must be lie in $A_a$, i.e., $2a=|x^2|=a$, and it follows that 
$a=0$. This contradiction shows that $x^2=0$. Similarly, we see that $y^2=0$.
Now consider the product $xy$ and we may assume that $xy=cx+dy$ for some $c,d\in k$. 
The left-hand side $xy$ is of degree $2a$ but the right-hand side has degree $a$, thus $cx+dy$ must be zero, i.e., $xy=0$.
A similar argument applies to $yx$ and obtains $yx=0$. Hence, this case produces a trivial left-symmetric color algebra. 

\textbf{Case 3}.  Suppose that $x\in A_0$ and $y\in A_a$ for some $a\neq 0$.  Note that
$x^2\in A_0$, $xy\in A_a$, and $yx\in A_a$. Moreover, $y^2\in A_aA_a\subseteq A_{2a}=A_0$, this is because
if $y^2\neq 0$, then it must be in either $A_0$ or $A_a$, however, as $2a\neq a$, we see that $y^2\notin A_a$.
Thus $y^2\in A_0$. Therefore, determining all 2-dimensional left-symmetric color algebras is actually equivalent to
determining all 2-dimensional left-symmetric superalgebras. The latter has been completed in \cite[Section 3]{ZB12}. 

\textbf{Case 4}.  Suppose that $x,y\in A_0$. Note that $A_0$ is a left-symmetric algebra and $\upepsilon(0,0)=1$, thus this situation is equivalent to classifying all 2-dimensional left-symmetric algebras; see for example
\cite[Section 5]{ZB12}.

\subsection{3-dimensional left-symmetric color algebras} 
A complete classification of 3-dimensional left-symmetric color algebras over a field (even the complex field) might be very complicated. We may see such complexity in the classifications of 
3-dimensional left-symmetric algebras \cite{Bai09} and 3-dimensional left-symmetric superalgebras \cite{ZB12}. Here, we don't intend to classify all 3-dimensional left-symmetric color algebras over a field but we may use the following approach to compute their cohomologies (actually, it should be working for any finite-dimensional left-symmetric color algebras). Let $A$ be a left-symmetric color algebras over $k$ and $V$ be an $A$-bimodule. 
\begin{enumerate}
  \item Consider the Lie color algebra $[A]$ and analyze the left $[A]$-module structure of $C^1(A,V)=\Hom(A,V)$;
  \item We may need to decompose $C^1(A,V)$ as a direct sum of some finite ``smaller" left $[A]$-modules, up to isomorphism;  or write $C^1(A,V)$ in a short exact sequences of left $[A]$-modules, and use the corresponding long exact sequence about cohomology spaces to compute $H^n([A], C^1(A,V))$; see \cite[Section 2]{SZ98};
  \item Use Theorem \ref{maint} to obtain $H^{n+1}(A,V)$.
\end{enumerate}

\begin{rem}{\rm
Note that in \cite{Sil97}, a classification (up to isomorphism) for 3-dimensional complex Lie color algebras with either 1-dimensional or zero homogeneous components has been reached, and using such classification, \cite{PS07} computed the cohomologies in trivial coefficients for 3-dimensional complex $\F_2^n$-graded Lie color algebras.  
\hbo}\end{rem}

We close this paper with the following example. 

\begin{exam}{\rm
Let's continue with the second Lie color algebra $L$ defined by the nonzero products
$$\{[x,y]=z,[z,x]=y,[y,z]=0\}$$
in Example \ref{exam2.10}. We would like to construct a left-symmetric color algebra $A$ such that $[A]=L$.
Note that $|x|=(1,1,0),|y|=(1,0,1)$, and $|z|=(0,1,1)$ and $L$ is spanned by $\{x,y,z\}$. We define $A$ to be given by
\begin{equation}
\label{ }
x^2=y^2=z^2=yx=xz=zy=yz=0, \quad  xy=z,\quad zx=y.
\end{equation}
A direct verification shows that $A$ is a left-symmetric color algebra and $[A]=L$. Suppose $V=\C$ denotes the trivial $A$-bimodule and we want to compute the cohomology spaces $H^i(A,\C)$ for
$i=0,1$. 

Note that $V=V_{(0,0,0)}$ is trivially graded by $\F_2^3$ and $\upepsilon(0,a)=1$ for all $a\in \F_2^3$. Thus it follows from 
Remark \ref{rem3.7} that $H^0(A,\C)=\C$.

We denote by $A^*:=C^1(A,\C)=\Hom(A,\C)$. By Example \ref{coh1}, $B^1(A,\C)\subseteq C^1(A,\C)=A^*$ consists of 
all homogeneous linear functions $\ell$ for which there exists some $c\in \C$ such that
$$\ell(r)=cr-\upepsilon(|c|,|r|)rc=cr-\upepsilon(0,|r|)rc=cr-rc=0$$ for  all $r\in A$.
Thus $B^1(A,\C)=0$. The cocycle space $Z^1(A,\C)\subseteq C^1(A,\C)=A^*$ consists of 
all homogeneous linear functions $\ell$ satisfying 
\begin{equation}
\label{ }
\ell(x_1x_2)=\ell(x_1)x_2+\upepsilon\left(|\ell|,|x_1|\right)x_1\ell(x_2),
\end{equation}
for all $x_1,x_2\in A$.  We write $\{x^*,y^*,z^*\}$ for a basis for $A^*$ dual to $\{x,y,z\}$. Taking $\ell$ to be any one of $\{x^*,y^*,z^*\}$, we observe that
\begin{equation}
\label{ }
\ell(xy)\neq \ell(x)y+\upepsilon\left(|\ell|,|x|\right)x\ell(y).
\end{equation}
Thus all $x^*,y^*,z^*$ are not contained in $Z^1(A,\C)$ and so $Z^1(A,\C)=0.$
Therefore, $H^1(A,\C)=0.$
\hbo}\end{exam}

\vspace{2mm}
\noindent \textbf{Acknowledgements}. The authors would like thank Professor Simon Xu for his helpful conversations and encouragement. This research was partially supported by the University of Saskatchewan
under grant No. APEF-121159.


\raggedright
\end{document}